\documentclass[10pt,journal]{IEEEtran}
\usepackage[utf8]{inputenc}
\usepackage{amsthm}
\usepackage[hidelinks]{hyperref}
\usepackage{extramath}
\usepackage[linesnumbered,ruled]{algorithm2e}
\usepackage{graphicx}
\usepackage{booktabs}
\usepackage{siunitx}
\usepackage[noadjust]{cite}
\usepackage{tikz}
\usepackage{tikztensor}
\usetikzlibrary{positioning}
\usetikzlibrary{patterns}
\usetikzlibrary{plotmarks}
\usetikzlibrary{shapes}
\usepackage{extramath}
\usepackage{ericmath}
\usepackage{url}
\usepackage{empheq}
\usepackage{pgfplots}
\pgfplotsset{compat=1.16}
\newtheorem{theorem}{Theorem}

\definecolor{front}{RGB}{114,199,254}
\definecolor{side}{RGB}{40,128,185}
\definecolor{top}{RGB}{76,167,232}
\definecolor{frontred}{RGB}{254,129,114}
\definecolor{sidered}{RGB}{185,88,40}
\definecolor{topred}{RGB}{232,107,76}

\makeatletter
\def\ps@IEEEtitlepagestyle{%
  \def\@oddfoot{\mycopyrightnotice}%
  \def\@oddhead{\hbox{}\@IEEEheaderstyle\leftmark\hfil\thepage}\relax
  \def\@evenhead{\@IEEEheaderstyle\thepage\hfil\leftmark\hbox{}}\relax
  \def\@evenfoot{}%
}
\def\mycopyrightnotice{%
  \begin{minipage}{\textwidth}
  \centering \scriptsize
  Copyright~\copyright~2022 IEEE. Personal use of this material is permitted. Permission from IEEE must be obtained for all other uses, in any current or future media, including\\reprinting/republishing this material for advertising or promotional purposes, creating new collective works, for resale or redistribution to servers or lists, or reuse of any copyrighted component of this work in other works by sending a request to pubs-permissions@ieee.org.
  \end{minipage}
}
\makeatother

\title{Canonical Polyadic Decomposition via the generalized Schur decomposition}

\author{
	 Eric~Evert$^{*\dagger}$, Michiel~Vandecappelle$^{*\dagger}$,~\IEEEmembership{Student member,~IEEE,} and Lieven~De~Lathauwer$^{*\dagger}$,~\IEEEmembership{Fellow,~IEEE}
	\thanks{\textbf{Funding:} Research supported by: (1) Flemish Government: Fonds de la Recherche Scientifique--FNRS and the Fonds Wetenschappelijk Onderzoek--Vlaanderen under EOS project no 30468160 (SeLMA) and Artificiële intelligentie (AI) Vlaanderen (3E190661); (2) KU Leuven Internal Funds C16/15/059 and ID-N project no 3E190402; (3) Leuven Institute for Artificial Intelligence (Leuven.ai).}	
	\thanks{$^*$KU Leuven, Dept. of Electrical Engineering ESAT/STADIUS, Kasteelpark Arenberg 10, bus 2446, B-3001 Leuven, Belgium}
	\thanks{$^\dagger$  KU Leuven -- Kulak, Group Science, Engineering and Technology, E. Sabbelaan 53, B-8500 Kortrijk, Belgium}
	\thanks{(Eric.Evert, Michiel.Vandecappelle, Lieven.DeLathauwer)@kuleuven.be.}
}
\begin{document}
	\IEEEtitleabstractindextext{%
		\begin{abstract}
		The canonical polyadic decomposition (CPD) is a fundamental tensor decomposition which expresses a tensor as a sum of rank one tensors. In stark contrast to the matrix case, with light assumptions, the CPD of a low rank tensor is (essentially) unique. The essential uniqueness of CPD makes this decomposition a powerful tool in many applications as it allows for extraction of component information from a signal of interest.

		One popular algorithm for algebraic computation of a CPD is the generalized eigenvalue decomposition (GEVD) which selects a matrix subpencil of a tensor, then computes the generalized eigenvectors of the pencil.  In this article, we present a simplification of GEVD which improves the accuracy of the algorithm. Surprisingly, the generalized eigenvector computation in GEVD is in fact unnecessary and can be replaced by a QZ decomposition which factors a pair of matrices as a product of unitary and upper triangular matrices. Computing a QZ decomposition is a standard first step when computing generalized eigenvectors, so our algorithm can been seen as a direct simplification of GEVD.
		\end{abstract}
		}
	
	\maketitle
	\IEEEdisplaynontitleabstractindextext
	
	\section{Introduction}
	\label{sec:intro} 
	
	Tensors, or multiindexed numerical arrays, are higher order generalizations of matrices and are natural structures for expressing data and signals which have inherent higher order structure. In this article we study the canonical polyadic decomposition (CPD) which expresses a tensor as a sum of rank one compoents. The CPD plays an important role in many applications due to the fact that, with mild assumptions, a low rank tensor has a unique CPD \cite{DD13,K77,SB00}. As such, one can recover underlying component information by computing a CPD of a low rank signal tensor \cite{Setal17,Cetal15,CJ10}. The essential uniqueness of CPD has helped make tensors and tensor methods common place in machine learning and signal processing \cite{Cetal15,Setal17}.
	
	A standard approach for computing a CPD of a low rank tensor is to first algebraically approximate the decomposition, then to refine the approximation with optimization routines. These algebraic approximations play an important role in optimization routines, as they are relatively inexpensive to compute and the strong initializations they provide can both improve final accuracy and reduce total computation time.  Notably, computing a best low rank approximation of a noisy low rank signal tensor is nonconvex and NP-hard \cite{HL13}, so good algebraic initializations greatly aid in getting reliable solutions.

A common approach for algebraic CPD computation is to first compute the generalized eigenvectors of a matrix subpencil of a tensor \cite{DD14,DD17,SK90,LRA93}. One of the tensor's factors can then be obtained by computing the inverse transpose of the matrix of generalized eigenvectors for the pencil. This generalized eigenvalue decomposition (GEVD) has been examined by many authors, e.g., see \cite{SK90,LRA93,FBK94}. 
	
	We present a simplification of the GEVD algorithm which is more accurate than the original algorithm. As it turns out, the generalized eigenvector computation in GEVD is not needed. Instead, one need only compute a generalized Schur decomposition of a subpencil of the tensor. This decomposition is also called the $\bQ \bZ$ decomposition. For a generic low rank tensor, the computed $\bQ$ and $\bZ$ will simultaneously upper triangularize all frontal slices of the tensor. Our key observation is that, once one has upper triangularized all frontal slices of the tensor, one can simply read off one of the tensors factor matrices from the diagonal entries of the upper triangular slices. Borrowing the perspective of \cite{EvertD19}, upper triangularizing the frontal slices of a low rank tensor reveals its ``joint generalized eigenvalues."
	
	Computing a $\bQ \bZ$ decomposition is a standard first step in computing generalized eigenvectors \cite{GV96}, so our algorithm can be seen as a direct simplification of GEVD. Intuitively, we view CPD computation as a (joint) generalized eigenvalue computation rather than a (joint) generalized eigenvector computation. The increase in accuracy of in our algorithm is due to the fact that computing the $\bQ \bZ$ decomposition only relies on unitary matrices as opposed to the general invertible matrices needed for generalized eigenvectors. Furthermore, the $\bQ \bZ$ method eliminates an inverse computation needed by GEVD.

It is worth noting that ours is not the first algorithm for CPD computation based on simultaneous upper triangularization of the frontal slices of a tensor. For example, in \cite{DDV04} frontal slices are jointly upper triangularized by minimizing the Frobenius norm of the lower triangular portion of the slices. However, the method used in \cite{DDV04} to obtain factor matrices after simultaneous upper triangularization is more involved.

Our $\bQ \bZ$ based CPD algorithm is presented in Section \ref{sec:qz}. In Section \ref{sec:lemmas} we prove that our algorithm successfully computes the CPD of a generic low rank tensor. The article ends with Section \ref{sec:exp} where we illustrate the performance of our algorithm on direction-of-arrival retrieval and in a fluorescence data experiment using the amino acid data set from \cite{bro1998multi}.

{\em Notation and terminology.} Let $\K$ denote either $\R$ or $\C$. We denote scalars, vectors, matrices, and tensors by lower case $(a)$, bold lower case $(\vec{a})$, bold upper case $(\bA)$, and calligraphic script $(\cA)$, respectively. For a matrix $\bA \in \K^{I \times J}$, we let $\bA^{\teT}$ denote the transpose of $\bA$ while $\bA^\teH$ denotes the conjugate transpose of $\bA$. If $\bA$ is invertible, then we let $\bA^{-\teT}$ denote the inverse of $\bA^\teT$. We use $D_r (\bA)$ denote the diagonal matrix whose diagonal entries are given by the $r$th row of $\bA$. We say a matrix is \df{generic} if it lies in a full measure set. Roughly speaking, a matrix is generic with probability equal to one. 

A \df{tensor} is a multiindexed array with entries in $\K$. The \df{order} of a tensor is the number of indices. Given a collection nonzero vectors $\bu^{(1)} \in \K^{I_1},\dots,\bu^{(N)} \in \K^{I_N}$, let 
\[
\bu^{(1)} \op \cdots \op \bu^{(N)}  \in \K^{I_1 \times \cdots \times I_N}
\]
denote the $I_1 \times \cdots \times I_N$ tensor with $i_1,i_2,\dots,i_N$ entry equal to $u^{(1)}_{i_1} u^{(2)}_{i_2} \cdots u^{(N)}_{i_N}$. A tensor of this form is called a \df{rank one tensor}. The minimal integer $R$ such that
\[
\mathcal{T}= \sum_{r=1}^R \bu^{(1)}_r \op \cdots \op \bu^{(N)}_r 
\]
where each $\bu_r^{(n)}$ has entries in $\K$ is called the $\K$-rank of the tensor $\mathcal{T}$, and a decomposition of this form is called a \df{canonical polyadic decomposition} (CPD) of $\mathcal{T}$. Compactly we write
$
\mathcal{T}=\cpdN. 
$
Here the matrix $\bU^{(n)} \in \K^{I_n \times R}$ has $\bu^{(n)}_r$ as its $r$th column. The matrix $\bU^{(n)}$ is called a \df{factor matrix} of $\mathcal{T}$.

A \df{mode-$\ell$ fiber} of a tensor is a vector obtained by fixing all indices but the $\ell$th. Defined for a third order tensor, the \df{mode-$\ell$ unfolding} $\mathbf{T}_{[\ell;j,i]}$ of $\mathcal{T}$ is the matrix obtained by stacking all mode-$\ell$ fibers of $\mathcal{T}$ as columns of a matrix, where the mode $i$ indices increment faster than the mode-$j$ indices. The \df{$\ell$-mode} product $\mathcal{T} \cdot_\ell \bA$ between a matrix $\bA$ and a tensor $\mathcal{T}$ is the tensor with mode-$\ell$ unfolding equal to $\bA \mathbf{T}_{[\ell;j,i]}$. We also make use of more general unfoldings for tensors of order greater than three. E.g., if $\mathcal{T} \in \K^{I_1 \times I_2 \times I_3 \times I_4}$, then $\mathbf{T}_{[4,3;2,1]}$ is an $I_1 I_2 \times I_3 I_4$ matrix with $I_1(i_2-1)+i_1, I_3(i_4-1)+i_3$ entry equal to $t_{i_1,i_2,i_3,i_4}$.

We often consider order three subtensors of a tensor. Given a tensor $\mathcal{T} \in \K^{I_1 \times \cdots \times I_N}$ and an integer $3 \leq n \leq N$ define
\[
\mathcal{T}[n] := \mathcal{T}(:,:,1,\dots,1,:,1,\dots,1)
\] 
where the third $:$ occurs in the $n$th mode of $\mathcal{T}$. That is, $\mathcal{T}[n]$ is an order three subtensor of $\mathcal{T}$ formed by fixing all but the first, second, and $n$th. In the case that $\mathcal{T}$ has order $3$, we call the matrices $\{\mathbf{t}(:,:,k)\}_{k=1}^{I_3}$ the \df{frontal slices} of $\mathcal{T}$.

	\section{CPD by QZ}
	\label{sec:qz}
	
	We now present the CPDQZ algorithm. We let $\cM \in \K^{I_1 \times \dots \times I_N}$ denote a measured tensor of interest, and we assume that $\cM$ has the form $\cM = \cT+\cN$ where $\cT$ is the signal portion of $\cM$ and where $\cN$ is noise. A standard assumption for generalized eigenvalue based algorithms such as GEVD and the generalized eigenspace decomposition (GESD), see \cite{EvertVD20}, is that the tensor $\mathcal{T}$ has at least two factor matrices with full column rank. In particular, letting $R$ denote the rank $\mathcal{T}$, one has $R \leq \min\{I_1,I_2\}$ up to a permutation of indices. Thus, by computing a (truncated) orthogonal compression of $\mathcal{\cM}$ such as a multilinear singular value decomposition, see e.g. \cite{DDV00b}, we can restrict to the case where $\mathcal{\cM}$ has rank $R$ and size $R \times R \times R_3 \times \cdots \times R_N$ with $R \geq R_n$ for all $n$.\footnote{We assume $R \geq R_n$ since tensor rank upper bounds multilinear rank.}
	
When compared to GEVD, our algorithm requires one additional mild assumption. Namely, we assume that the matrix $\bU^{(N)} \odot \dots \odot \bU^{(3)}$ has full column rank. Here $\odot$ denotes the Khatri-Rao product. In the upcoming CPDQZS variation of our algorithm, one instead must make the stronger assumption that there is an index $3 \leq n \leq N$ such that $R_n=R$ and $\bU^{(n)}$ is invertible. These deterministic assumptions are satisfied for generic factor matrices provided the matrices in question all have at least as many rows as columns (before orthogonal compression).

The key observation behind CPDQZ is that one may compute all but two factor matrices of a tensor using a single $\bQ \bZ$ decomposition together with $n$-mode products. In particular, given a rank $R$ tensor $\mathcal{T} \in \K^{R \times R \times R_3 \times \cdots \times R_N} $ which satisfies our assumptions, to obtain factor matrices $\bU^{(n)}$ for $n \geq 3$, one need only compute unitary matrices $\bQ$ and $\bZ$ such that 
\[
\bQ \mathbf{T}[3](:,:,1) \bZ  \qquad  \bQ \mathbf{T}[3](:,:,2) \bZ
\]
is a $\bQ \bZ$ decomposition of the matrix pencil\footnote{Linear combinations of matrices $\mathbf{T}(:,:,i_3,\dots,i_N)$ may be taken to form the matrix pencil used for $\bQ \bZ$ computation. We use a fixed pencil for the sake of exposition. A popular pencil choice which can improve numerics is the first two slices of the (truncated) core of the MLSVD of $\cT$. However, this does not solve all issues, e.g., see \cite[Section 4.4]{EvertVD20}.}
\[
(\mathbf{T}[3](:,:,1),\mathbf{T}([3](:,:,2)),
\]
i.e., such that  these matrices are both upper triangular.\footnote{If a matrix pencil $(\bM_1,\bM_2)$ has complex generalized eigenvalues and $\bQ$ and $\bZ$ have real entries, then there will be $2 \times 2$ blocks on the diagonal of each matrix $\bQ \bM_i \bZ$. However, if $\bM_1$ and $\bM_2$ are formed from linear combinations of frontal slices of a real rank $R$ tensor $\cT$ which meets our assumptions, then all generalized eigenvalues of $(\bM_1,\bM_2)$ are real, and each $\bQ \bM_i \bZ$ will be upper triangular, e.g., see \cite{LRA93}. Thus, the presence of $2 \times 2$ blocks on the diagonal may indicate that the low real rank model chosen is not appropriate or that the signal to noise ratio is not sufficiently high.}

Define $\mathcal{T}_{qz}:= \mathcal{T} \cdot_1 \bQ \cdot_2 \bZ^\teT$. It is then a matter of technical formula manipulation to show that
\[
\bU^{(n)} (r,:) = \mathbf{t}_{qz}[n](r,r,:)
\]
for all $n = 3,\dots, N$ and all $r = 1,\dots, R.$ See the supplementary materials for details.

It remains to compute $\bU^{(1)}$ and $\bU^{(2)}$. To do this we first compute $\bU^{(2)} \odot \bU^{(1)}$ by solving the overdetermined system
\[
\mathbf{T}_{[N,N-1,\dots,3;2,1]} = (\bU^{(N)} \odot \dots \odot \bU^{(3)}) (\bU^{(2)} \odot \bU^{(1)}) ^\teT.
\]
The $r$th column of the matrix $\bU^{(2)} \odot \bU^{(1)}$ is the Kronecker product of the $r$th column of $\bU^{(2)}$ with that of $\bU^{(1)}$ and can therefore be reshaped into a rank-$1$ matrix. It follows that the columns of these factors can be obtained by computing a rank-$1$ approximation of each reshaped column of $\bU^{(2)} \odot \bU^{(1)}$.

As an alternative, suppose $\bU^{(n)} \in \K^{R \times R}$ has full rank for some $n \in \{3,\dots,N\}$. For ease of exposition we take $n=N$. Then the remaining factor matrices can be computed by first solving for $(\bU^{(N-1)}\odot \dots \odot \bU^{(1)})$ in
\[
\mathbf{T}_{[N-1,\dots,1;N]} = (\bU^{(N-1)}\odot \dots \odot \bU^{(1)}) (\bU^{(N)})^\teT,
\]
then computing rank-$1$ approximations of the appropriately reshaped columns of this matrix.  We call the first approach CPDQZ and the second CPDQZS(ingle). For emphasis, CPDQZ and CPDQZS are the same for order three tensors. 

For tensors or order greater than three, the difference between these methods is that the CPDQZS method uses rank-$1$ tensor approximations to compute most of its factors, while the CPDQZ relies more heavily on the initial $\bQ \bZ$ decomposition to obtain factors. Since a best rank-$1$ tensor approximation can often be accurately computed and since there is less opportunity for error accumulation before rank-$1$ tensor approximations are employed in CPDQZS, the CPDQZS method is expected to be more accurate than CPDQZ. This expectation is supported by our numerical experiments.

The growth rate of the cost of both GEVD and CPDQZ is $\mathcal{O}(R^3 (\Pi_{i=3}^{N} R_i))$. The growth rate of the cost of CPDQZS is $\mathcal{O}(R^4 (\Pi_{i=3}^{N-1} R_i))$. If $R_N=R$, then these growth rates coincide. However, the coefficients of the costs for the methods can be very different. The algorithms all share a step whose cost grows at $\mathcal{O}(R^4 (\Pi_{i=3}^{N-1} R_i))$. CPDQZS and GEVD each have additional steps with this growth rate, while CPDQZ does not. Thus, the coefficient of the cost of CPDQZ is lower than that of GEVD or CPDQZS. In practice, CPDQZ is observed to be much faster than CPDQZS and GEVD, see Section \ref{sec:exp}. See the supplementary materials for further discussion.

\section{Algorithm derivation}
\label{sec:lemmas}

We now derive the CPDQZ algorithm. To ease exposition, we temporarily assume $\cT \in \K^{R \times R \times R}$ has order three. The first observation used in the derivation of the CPDQZ algorithm is that in the special case where $\cT = [\![\bU^{(1)},\bU^{(2)},\bU^{(3)}]\!]$ has upper triangular frontal slices and meets our assumptions, then in an appropriate column ordering, the matrix $\bU^{(1)}$ is upper triangular while the matrix $\bU^{(2)}$ is lower triangular.  

This fact follows quickly from the main observation behind GEVD. Namely, if $\bU^{(1)}$ and $\bU^{(2)}$ are $R \times R$ invertible matrices, then the columns of the matrix $(\bU^{(2)})^{-\teT}$ are equal to the generalized eigenvectors of the matrix pencil 
\[
(\bT(:,:,1),\bT(:,:,2)).
\]
In the case that $\bT(:,:,1)$ and $\bT(:,:,2)$ are upper triangular, a routine argument shows that the pencil's generalized eigenvectors can be ordered so that $(\bU^{(2)})^{-\teT}$ is upper triangular, hence $\bU^{(2)}$ is lower triangular. Having shown that $\bU^{(2)}$ is lower triangular, one may use the assumption that $\cT$ has upper triangular frontal slices together with the formula
\beq
\label{eq:FrontalSliceFormula}
\bT(:,:,r) = \bU^{(1)} D_r (\bU^{(3)}) (\bU^{(2)})^\teT \qquad  \mathrm{for  \ } r=1,\dots,R
\eeq
to conclude not only that $\bU^{(1)}$ is upper triangular, but also that, up to scaling, one has
\[
\bU^{(3)} (:,r) = \text{diag}(\mathbf{T}(:,:,r)) \qquad  \mathrm{for  \ } r=1,\dots,R.
\]

From this point we need only show that if $\cT  = [\![ \bU^{(1)},\bU^{(2)},\bU^{(3)}]\!]\in \K^{R \times R \times R}$ is an arbitrary tensor which meets our assumptions and if $\bQ$ and $\bZ$ are matrices which give a $\bQ \bZ$ decomposition of the matrix pencil 
\[
(\bQ \bT(:,:,1) \bZ,\bQ \bT(:,:,2) \bZ),
\]
then the tensor $\cT \cdot_1 \bQ \cdot_2 \bZ^\teT = [\![ \bQ \bU^{(1)},\bZ^\teT \bU^{(2)},\bU^{(3)}]\!]$ has upper triangular frontal slices. As we shall explain, this fact follows from repeating preceding argument. 

Our assumptions guarantee that the pencil $(\bT(:,:,1),\bT(:,:,2))$ has real generalized eigenvalues \cite{LRA93}, hence $\bQ \bT(:,:,1)\bZ$ and $\bQ \bT(:,:,2) \bZ$ are upper triangular by definition of the $\bQ \bZ$ factorization and the matrix of generalized eigenvectors of this pencil, i.e, the inverse transpose of the second factor matrix of the tensor $\bT\cdot_1 \bQ \cdot_2 \bZ^{\teT}$,  can be taken to be upper triangular, e.g., see \cite{GV96}.  That is, $\bZ^\teT \bU^{(2)}$ can be taken to be lower triangular. Applying equation \eqref{eq:FrontalSliceFormula} to $\bT\cdot_1 \bQ \cdot_2 \bZ^{\teT}$ shows that $\bQ \bU^{(1)}$ is upper triangular and that all frontal slices of this tensor are upper triangular. Since $\cT$ and $\cT \cdot_1 \bQ \cdot_2 \bZ^\teT$ have the same third factor matrix it then easily follows that the CPDQZ algorithm successfully recovers the CPD of a tensor which meets the assumptions stated at the beginning of Section \ref{sec:qz}.

The extension of CPDQZ and CPDQZS to tensors of order greater than three follows a routine argument using the fact that an order $N$ tensor with CPD $\cpdN$ can be reshaped to an order three tensor with CPD $[\![ \bU^{(1)},\bU^{(2)},\bU^{(N)}  \odot \dots \odot \bU^{(3)}]\!]$.

	\section{Experiments}
	\label{sec:exp}
	The proposed QZ methods make different trade-offs with respect to accuracy and speed. In this section, we compare both methods with the classical GEVD algorithm and the more recent GESD algorithm~\cite{evert2020recursive}. We use a machine with an AMD Ryzen 5~5600H CPU at 3,30GHz and 16GB of RAM using MATLAB R2021b and Tensorlab 3.0~\cite{vervliet2016tensorlab}.
	
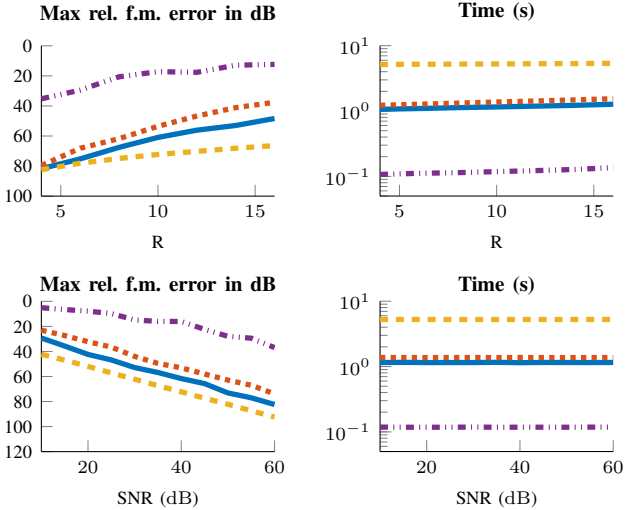
\begin{figure}
		\centering
%
%
\definecolor{mycolor1}{rgb}{0.00000,0.44700,0.74100}%
\definecolor{mycolor2}{rgb}{0.85000,0.32500,0.09800}%
\definecolor{mycolor3}{rgb}{0.92900,0.69400,0.12500}%
\definecolor{mycolor4}{rgb}{0.49400,0.18400,0.55600}%
\pgfplotsset{every tick label/.append style={font=\scriptsize}}
\begin{tikzpicture}

\begin{axis}[%
width=3.1cm,
height=2.0cm,
at={(0cm,3.4cm)},
scale only axis,
xmin=4,
xmax=16,
xlabel={\scriptsize R},
ymode=log,
ymin=1e-05,
ymax=1e-0,
ytick = {1e-0,1e-1,1e-2,1e-3,1e-4,1e-5},
yticklabels = {0,20,40,60,80,100},
yminorticks=true,
title={Max rel.\ f.m.\ error in dB},
axis background/.style={fill=white},
title style={font=\bfseries\footnotesize},
axis x line*=bottom,
axis y line*=left,
legend style={legend cell align=left,align=left,draw=white!15!black}
]
\addplot [color=mycolor1,solid,line width = 2pt]
  table[row sep=crcr]{%
4	8.17442774519231e-05\\
6	0.000173545659480016\\
8	0.000418982563865303\\
10	0.000897377501204146\\
12	0.00156234525408739\\
14	0.00223650627041485\\
16	0.00382835881749107\\
};
\label{CPDQZS}
\addplot [color=mycolor4,dash dot dot,line width = 2pt]
  table[row sep=crcr]{%
4	0.0172456040458189\\
6	0.0337150791967955\\
8	0.093192581659901\\
10	0.137995299798376\\
12	0.131121848988428\\
14	0.226221355287014\\
16	0.240926004899586\\
};
\label{CPDQZ}
\addplot [color=mycolor2,dotted,line width = 2pt]
  table[row sep=crcr]{%
4	0.000104332324421058\\
6	0.000391942932508538\\
8	0.000817279566359919\\
10	0.00210707077455458\\
12	0.00461505929317431\\
14	0.00878825730125472\\
16	0.0131840987860146\\
};
\label{gevd}
\addplot [color=mycolor3,dashed,line width = 2pt]
  table[row sep=crcr]{%
4	7.5420760144031e-05\\
6	0.000124256721329953\\
8	0.000179393910424867\\
10	0.000243624001106932\\
12	0.000308604132407451\\
14	0.000391461934538794\\
16	0.000478556821444878\\
};
\label{gesd}
\end{axis}

\begin{axis}[%
width=3.1cm,
height=2.0cm,
at={(4.5cm,3.4cm)},
scale only axis,
xmin=4,
xmax=16,
xlabel={\scriptsize R},
ymode=log,
ymin=0.05,
ymax=10,
yminorticks=true,
title={Time (s)},
axis background/.style={fill=white},
title style={font=\bfseries\footnotesize},
axis x line*=bottom,
axis y line*=left,
legend style={legend cell align=left,align=left,draw=white!15!black}
]
\addplot [color=mycolor1,solid,line width = 2pt]
  table[row sep=crcr]{%
4	1.0630163\\
6	1.0922802\\
8	1.12694825\\
10	1.15508895\\
12	1.18836775\\
14	1.22120865\\
16	1.2710733\\
};

\addplot [color=mycolor4,dash dot dot,line width = 2pt]
  table[row sep=crcr]{%
4	0.10712875\\
6	0.11134435\\
8	0.1139644\\
10	0.11814915\\
12	0.12252305\\
14	0.12748435\\
16	0.13626295\\
};

\addplot [color=mycolor2,dotted,line width = 2pt]
  table[row sep=crcr]{%
4	1.2257713\\
6	1.27217435\\
8	1.32134035\\
10	1.3762582\\
12	1.42530775\\
14	1.48168315\\
16	1.54154535\\
};

\addplot [color=mycolor3,dashed,line width = 2pt]
  table[row sep=crcr]{%
4	5.15349275\\
6	5.2004556\\
8	5.22583305\\
10	5.2616718\\
12	5.2970096\\
14	5.3265306\\
16	5.3708199\\
};

\end{axis}

\begin{axis}[%
	width=3.1cm,
	height=2.0cm,
	at={(0cm,0cm)},
	scale only axis,
	xmin=10,
	xmax=60,
	xlabel={\scriptsize SNR (\si{\decibel})},
	ymode=log,
	ymin=1e-6,
	ytick = {1e-0,1e-1,1e-2,1e-3,1e-4,1e-5,1e-6},
	yticklabels = {0,20,40,60,80,100,120},
	ymax=1,
	yminorticks=true,
	title={Max rel.\ f.m.\ error in dB},
	axis background/.style={fill=white},
	title style={font=\bfseries\footnotesize, at={(0.5,0.9)}},
	axis x line*=bottom,
	axis y line*=left,
	legend style={legend cell align=left,align=left,draw=white!15!black}
	]
	\addplot [color=mycolor1,solid,line width = 2pt]
	table[row sep=crcr]{%
	10	0.0343498461865444\\
15	0.0163723769857914\\
20	0.00758824689851713\\
25	0.00457516230884916\\
30	0.00226831994471322\\
35	0.00144168813290182\\
40	0.000825119531856585\\
45	0.000519128843482251\\
50	0.000220308933937591\\
55	0.000141063003611483\\
60	7.60694599361006e-05\\
	};
	\addplot [color=mycolor4,dash dot dot,line width = 2pt]
	table[row sep=crcr]{%
	10	0.555038762470017\\
15	0.469738069056889\\
20	0.404198204546676\\
25	0.325515557317449\\
30	0.180408128180455\\
35	0.157363261991236\\
40	0.157612376857481\\
45	0.0757538425986829\\
50	0.0400238334954679\\
55	0.0336711244712485\\
60	0.014240032301295\\
	};
	\addplot [color=mycolor2,dotted,line width = 2pt]
	table[row sep=crcr]{%
		10	0.0704735918867488\\
15	0.0428957337981458\\
20	0.024402881131278\\
25	0.0148571151196085\\
30	0.00632663793507694\\
35	0.00333760275287516\\
40	0.00220611557557422\\
45	0.00125300078358425\\
50	0.000715906357667629\\
55	0.000444690481396178\\
60	0.000202579668944919\\
	};
	\addplot [color=mycolor3,dashed,line width = 2pt]
	table[row sep=crcr]{%
	10	0.00768714778776904\\
15	0.00448050755001364\\
20	0.0025122156717023\\
25	0.00138552081182312\\
30	0.000777612328522113\\
35	0.000427620157387879\\
40	0.000243934035374067\\
45	0.000136031256723392\\
50	7.77635498685349e-05\\
55	4.26345502927623e-05\\
60	2.40046654478079e-05\\
	};
\end{axis}

\begin{axis}[%
	width=3.1cm,
	height=2.0cm,
	at={(4.5cm,0cm)},
	scale only axis,
	xmin=10,
	xmax=60,
	xlabel={\scriptsize SNR (\si{\decibel})},
	ymode=log,
	ymin=0.05,
	ymax=10,
	yminorticks=true,
	title={Time (s)},
	axis background/.style={fill=white},
	title style={font=\bfseries\footnotesize, at={(0.5,0.88)}},
	axis x line*=bottom,
	axis y line*=left,
	legend style={legend cell align=left,align=left,draw=white!15!black}
	]
	\addplot [color=mycolor1,solid,line width = 2pt]
	table[row sep=crcr]{%
	10	1.1500502\\
15	1.15901335\\
20	1.15022115\\
25	1.1521382\\
30	1.1508732\\
35	1.16024505\\
40	1.14564155\\
45	1.15553675\\
50	1.1527224\\
55	1.15149125\\
60	1.1548991\\
	};
	
	\addplot [color=mycolor4,dash dot dot,line width = 2pt]
	table[row sep=crcr]{%
	10	0.11823105\\
15	0.11765175\\
20	0.1177149\\
25	0.11770295\\
30	0.11840705\\
35	0.11766695\\
40	0.11784255\\
45	0.1177947\\
50	0.1176742\\
55	0.1182116\\
60	0.11802235\\
	};
	
	\addplot [color=mycolor2,dotted,line width = 2pt]
	table[row sep=crcr]{%
		10	1.3735746\\
15	1.37483825\\
20	1.3699639\\
25	1.37342435\\
30	1.37524305\\
35	1.37041105\\
40	1.3745698\\
45	1.3729006\\
50	1.37084445\\
55	1.3706539\\
60	1.37574765\\
	};
	
	\addplot [color=mycolor3,dashed,line width = 2pt]
	table[row sep=crcr]{%
10	5.2569238\\
15	5.2512\\
20	5.24854605\\
25	5.24546815\\
30	5.25173535\\
35	5.25012255\\
40	5.25111425\\
45	5.25372765\\
50	5.2484466\\
55	5.25272185\\
60	5.2535339\\
	};
	
\end{axis}
\end{tikzpicture}%
		\caption{\label{fig:varying}The QZ methods make different trade-offs with respect to time and accuracy. Over $50$ trials, CPDQZ (\ref{CPDQZ}) is the fastest, but least accurate, while GESD (\ref{gesd}) is the most accurate and the slowest. CPDQZS (\ref{CPDQZS}) and GEVD (\ref{gevd}) are equally fast, but the former is more accurate. The relative performance of the methods for the fourth-order tensors is consistent over varying ranks (with SNR $\SI{40}{\decibel}$) and SNRs (with rank $R=10$). In both cases, CPDQZS is on average about $8$ dB more accurate than GEVD.}
	\end{figure}

	In a first experiment, we generate fourth-order low-rank tensors $\mathcal{T}$ by sampling the entries of factor matrices $\mathbf{A}$, $\mathbf{B}$, $\mathbf{C}$ and $\mathbf{D}$ of dimensions $80 \times R$ from the uniform distribution on $[0,1]$ and we normalize all columns to unit length. The rank $R$ is varied in the range $[4,16]$ and Gaussian noise is added such that the SNR is $\SI{40}{\decibel}$. For the estimated factor matrices $\hat{\mathbf{A}}$, $\hat{\mathbf{B}}$, $\hat{\mathbf{C}}$ and $\hat{\mathbf{D}}$, we show maximal relative factor matrix errors compared to the true factors, defined as 
	\[
	\max\left(\frac{\Vert \mathbf{A} - \hat{\mathbf{A}} \Vert_\text{F}}{\Vert \mathbf{A}\Vert_\text{F}},\frac{\Vert \mathbf{B} - \hat{\mathbf{B}} \Vert_\text{F}}{\Vert \mathbf{B}\Vert_\text{F}},\frac{\Vert \mathbf{C} - \hat{\mathbf{C}} \Vert_\text{F}}{\Vert \mathbf{C}\Vert_\text{F}},\frac{\Vert \mathbf{D} - \hat{\mathbf{D}} \Vert_\text{F}}{\Vert \mathbf{D}\Vert_\text{F}}\right),
	\]
	where columns of $\hat{\mathbf{A}}$, $\hat{\mathbf{B}}$, $\hat{\mathbf{C}}$ and $\hat{\mathbf{D}}$ and have been optimally permuted and scaled to match the columns of $\mathbf{A}$, $\mathbf{B}$, $\mathbf{C}$ and $\mathbf{D}$, respectively. The results are shown in the top plots of Figure~\ref{fig:varying}. CPDQZ is the fastest method, but is less accurate compared to CPDQZS and GEVD. CPDQZS is as fast as GEVD, but is more accurate over the whole range of ranks. In the bottom plots of Figure~\ref{fig:varying}, the experiment is repeated, but now the SNR is varied between $10$ and $\SI{60}{\decibel}$, with the rank $R=10$ fixed. The same relative performance is seen for the four methods.
	
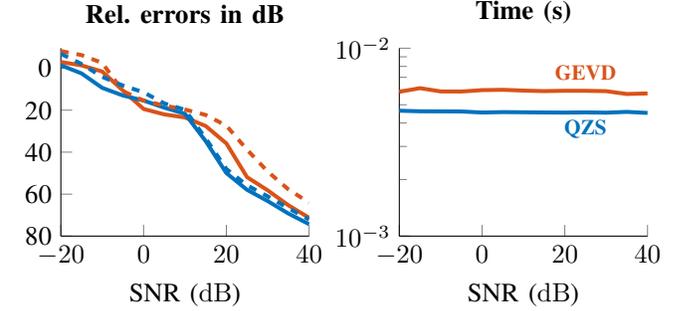
\begin{figure}
    	\centering
%
%
\definecolor{mycolor1}{rgb}{0.00000,0.44700,0.74100}%
\definecolor{mycolor2}{rgb}{0.85000,0.32500,0.09800}%
\definecolor{mycolor3}{rgb}{0.92900,0.69400,0.12500}%
\definecolor{mycolor4}{rgb}{0.49400,0.18400,0.55600}%
\begin{tikzpicture}

\begin{axis}[%
width=3.3cm,
height=2.5cm,
at={(0cm,0cm)},
scale only axis,
xmin=-20,
xmax=40,
xlabel={SNR ($\si{\decibel}$)},
ymode=log,
ymin=0.0001,
ymax=3,
ytick = {1e-0,1e-1,1e-2,1e-3,1e-4},
yticklabels = {0,20,40,60,80},
yminorticks=true,
axis background/.style={fill=white},
title style={font=\bfseries},
title={Rel. errors in dB},
axis x line*=bottom,
axis y line*=left,
legend style={legend cell align=left,align=left,draw=white!15!black}
]
\addplot [color=mycolor2,solid,line width=1.5pt]
  table[row sep=crcr]{%
-20	1.40712290190521\\
-15	1.16605310324248\\
-10	0.834345742393365\\
-5	0.293623934632326\\
0	0.106871794300969\\
5	0.0793449299129675\\
10	0.067034130266788\\
15	0.0420109516808151\\
20	0.0161585377181663\\
25	0.00253773771422795\\
30	0.00122899687734663\\
35	0.000542857000653001\\
40	0.000273328307258484\\
};
\label{gevdazi}
\addplot [color=mycolor1,solid,line width=1.5pt]
  table[row sep=crcr]{%
-20	1.17691240812226\\
-15	0.751653564772446\\
-10	0.340134242973551\\
-5	0.223988039307387\\
0	0.169802582721228\\
5	0.114892876975422\\
10	0.0819671167619221\\
15	0.0183475353279176\\
20	0.00317095153404353\\
25	0.00126695331733697\\
30	0.000682860255574686\\
35	0.000344506044564039\\
40	0.000192475802049853\\
};
\label{qzazi}
\addplot [color=mycolor2,dashed,line width=1.5pt]
  table[row sep=crcr]{%
-20	2.56156652579874\\
-15	2.05325983754165\\
-10	1.33759526631107\\
-5	0.283955356702785\\
0	0.169392671825325\\
5	0.128578349079894\\
10	0.103605037548347\\
15	0.0769213170538817\\
20	0.0428204868185539\\
25	0.0114492586183623\\
30	0.00344245936481209\\
35	0.00129629140748166\\
40	0.000627569619446287\\
};
\label{gevdele}
\addplot [color=mycolor1, dashed,line width=1.5pt]
  table[row sep=crcr]{%
-20	2.20629030488283\\
-15	1.26060175191087\\
-10	0.621572360191161\\
-5	0.386046498858233\\
0	0.266746571672884\\
5	0.144836544800713\\
10	0.102225406970788\\
15	0.021342235082441\\
20	0.00398143170669914\\
25	0.00162054810985345\\
30	0.000884085756771264\\
35	0.000460163525108553\\
40	0.000249789118037782\\
};
\label{qzele}
\end{axis}

\begin{axis}[%
width=3.3cm,
height=2.5cm,
at={(4.5cm,0cm)},
scale only axis,
xmin=-20,
xmax=40,
xlabel={SNR ($\si{\decibel}$)},
ymin=0.001,
ymax=0.01,
ymode = log,
axis background/.style={fill=white},
title style={font=\bfseries},
title={Time (s)},
axis x line*=bottom,
axis y line*=left,
legend style={legend cell align=left,align=left,draw=white!15!black}
]
\addplot [color=mycolor2,solid,line width=1.5pt]
  table[row sep=crcr]{%
-20	0.005863\\
-15	0.00612905\\
-10	0.00588805\\
-5	0.00588145\\
0	0.00598295\\
5	0.0060071\\
10	0.00594935\\
15	0.0059103\\
20	0.00593665\\
25	0.00593345\\
30	0.00590705\\
35	0.0057055\\
40	0.00574575\\
};

\addplot [color=mycolor1,solid,line width=1.5pt]
  table[row sep=crcr]{%
-20	0.00464875\\
-15	0.00461595\\
-10	0.00461125\\
-5	0.0046053\\
0	0.0045471\\
5	0.0045705\\
10	0.00456075\\
15	0.0045494\\
20	0.00454355\\
25	0.00455515\\
30	0.0045381\\
35	0.00458335\\
40	0.00452395\\
};
\node [mycolor2] at (25,0.00745){\footnotesize \bfseries GEVD};
\node [mycolor1] at (25,0.0037){\footnotesize \bfseries QZS};
\end{axis}
\end{tikzpicture}%
    	\caption{\label{fig:doa} CPDQZS estimates the elevations of the sources more accurately than GEVD outside of a small SNR window. In addition, CPDQZS has a lower computational cost. Left: mean relative errors over the eight sources of the azimuths (\ref{qzazi}) and elevations (\ref{qzele}) that are estimated for CPDQZS and the azimuths (\ref{gevdazi}) and elevations (\ref{gevdele}) that are estimated for GEVD. Right: computation time for both methods. Medians over $500$ trials are shown.}
    \end{figure}

In a more applied experiment, the QZ approach is compared to GEVD in a direction-of-arrival (DOA) retrieval experiment for line-of-sight signals impinging on a uniform rectangular array (URA). The CPD can be applied to find the DOAs in this case~\cite{sidiropoulos2000parallel,roy1989esprit,HRG08}. The URA has $M\times M$ sensors, where $M=20$, and collects $K=20$ samples from $R=8$ far-field sources. These have azimuths $\begin{smallmatrix}[15&20&25&30&35&40&50&55]\end{smallmatrix}$ and elevations $[\begin{smallmatrix}5&10&20&25&30&35&40&45]\end{smallmatrix}$ degrees, respectively. We consider omnidirectional sensors, evenly-spaced with inter-sensor spacing $\Delta$. The azimuths $\mathbf{Z}\in\mathbb{C}^{M\times R}$ and elevations $\mathbf{L}\in\mathbb{C}^{M\times R}$ of the $R$ sources that we collect at each sensor yield an observed tensor $\mathcal{T}\in\mathbb{C}^{M\times M \times K}$, where the $k$th frontal slice has a low rank decomposition  $\bT_k=\mathbf{A}^{(k)} \mathrm{diag}( \mathbf{s}^{(k)}) (\mathbf{E}^{(k)})^\teT$. The matrix $\mathbf{A}^{(k)}\in\mathbb{C}^{M\times R}$ has entries $a^{(k)}_{mr} = \exp({(m-1)2\pi/\lambda\sin(z_{rk}\pi/180)\Delta i})$, while $\mathbf{E}^{(k)}\in\mathbb{C}^{M\times R}$ has entries $e^{(k)}_{mr}=\exp({(m-1)2\pi/\lambda\sin(l_{rk}\pi/180)\Delta i})$. The vector $\mathbf{s}^{(k)}\in\mathbb{C}^{R}$ holds the sources and $\lambda$ is the signal wavelength. By computing a rank-$R$ CPD of a few frontal slices of $\mathcal{T}$, the azimuths and elevations of the sources can be recovered from the first two CPD factors. The tensor $\mathcal{T}$ is perturbed with Gaussian distributed noise in the range $[-20,40]\si{\decibel}$ and we compare the accuracy and speed of the CPDs obtained with GEVD and CPDQZS. In Figure~\ref{fig:doa}, one can find the median error (left) and computation time (right) over $500$ trials for the azimuth and elevation estimations of the eight sources. QZ is more accurate than GEVD, especially for the estimation of the source elevations, and is faster as well.

\begin{figure}
		\centering
		\input{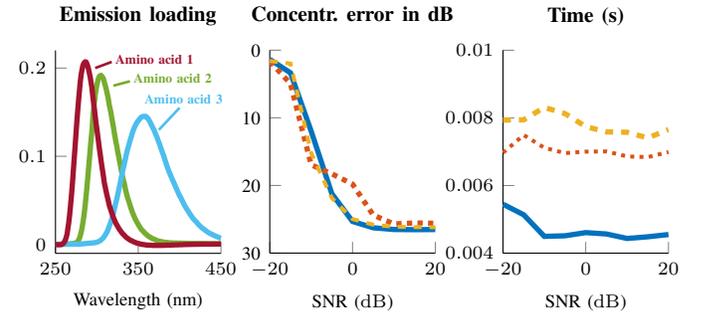}
		\caption{\label{fig:spectra}The accuracy of CPDQZS (\ref{CPDQZS}) for the estimation of the amino acid concentrations is as least as good as GEVD (\ref{gevd}) and GESD (\ref{gesd}), while requiring a lower computation time. The theoretical emission loading for the three amino acids in function of the wavelength is shown in the left plot. In the center plot, the mode-$1$ factor matrix error is shown, which corresponds to the estimated concentrations of the amino acids in the five mixtures. The right plot shows the computation times. All results are medians over $100$ trials.}
	\end{figure}

	Lastly, we use the amino acid dataset from \cite{bro1998multi}. This dataset holds the measured emission and excitation spectra of five amino acid mixtures in a tensor $\mathcal{T}$ of size $5\times 201 \times 61$. The first mode corresponds to the five mixtures and the second and third modes refer to the emission and excitation spectra, respectively. In the left plot of Figure~\ref{fig:spectra}, the theoretical emission spectra are shown for the three amino acids. Following \cite{bro1998multi}, $\mathcal{T}$ admits a rank-$3$ CPD with each component corresponding to one amino acid. We add Gaussian noise with SNRs between $-20$ and $\SI{20}{\decibel}$ and use GEVD, GESD and CPDQZS to compute a CPD of the perturbed $\mathcal{T}$. We compare the accuracy of the estimated concentrations of the amino acids in each mixture, provided by the mode-$1$ factor matrix of the CPD. The relative mode-$1$ factor matrix errors are shown in the center plot of Figure~\ref{fig:spectra}, while the computation times are shown in the right plot. CPDQZS is about as accurate as GEVD and GESD over the range of SNRs and is also notably faster.

	\section{Conclusion}
	\label{sec:conc}

	We presented the novel CPDQZ and CPDQZS algorithms for algebraic CPD computation which can be viewed as direct simplifications of the popular GEVD algorithm. These algorithms replace the generalized eigenvector computation of GEVD with a $\bQ \bZ$ decomposition which is used to upper triangularize the frontal slices of a tensor. One factor matrix is then obtained by reading the diagonal entries of the upper triangular slices. We showed in experiments that, in the case of CPDQZS, this simplification results in an increase in accuracy when compared to GEVD. For tensors of order four or more\footnote{ Recall that CPDQZ and CPDQZS are identical for tensors of order three.}, CPDQZ is observed to be less accurate but faster than GEVD.
	
	In future work we will investigate a deflation style algorithm in the spirit of GESD where the deflation step is based on a $\bQ \bZ$ computation rather than generalized eigenspace computations.

	\bibliographystyle{ieeetran}
	\bibliography{references,../References/references}
	
	
\end{document}